\documentclass[runningheads]{llncs}
\usepackage[T1]{fontenc}
\usepackage{graphicx}
\begin{document}
\title{ MathPartner is a breakthrough technology for natural sciences education, scientific and engineering applications}
\titlerunning{MathPartner is a breakthrough technology}
%
\author{Gennadi Malaschonok  \orcidID{0000-0002-9698-6374}, Roman Sakh \thanks{Preprint of the paper:\   Malaschonok, G., Sakh, R. (2025). MathPartner Is a Breakthrough Technology for Natural Sciences Education, Scientific and Engineering Applications. In: Ermolayev, V., et al. Information and Communication Technologies in Education, Research, and Industrial Applications. ICTERI 2024. Communications in Computer and Information Science, vol 2359. Springer, Cham. https://doi.org/10.1007/978-3-031-81372-6\_15} }
\authorrunning{Gennadi Malaschonok, Roman Sakh }
%
\institute{National University “Kyiv-Mohyla Academy”, Kyiv, Ukraine  
\email{malaschonok@ukma.edu.ua, roman.sakh@ukma.edu.ua}\\
\url{https://www.fin.ukma.edu.ua/department-network-technologies}  \\
 }
\maketitle              
\begin{abstract}
The article provides a brief description of the MathPartner service. This freely available cloud-based Mathematics is a universal system for symbolic-numeric calculations. Its Mathpar language is a subset of the LaTeX language, but allows you to create mathematical texts that contain “computable” mathematical operators.
This opens up completely new opportunities for improving the educational process for all natural science disciplines, for the use of mathematics in scientific and engineering calculations.
To save and freely exchange educational and other texts in the Mathpar language, a GitHub repository has been created.
 It is concluded that cloud mathematics MathPartner is a new breakthrough technology for school and university natural science education, for scientific and engineering applications.

\end{abstract}

%
%
\section{Introduction}
 
 Today, traditional means of learning mathematics, science and engineering are widely used. Educational tasks are solved with the help of a pen and notebook. The effectiveness of the educational process remains low.

It can be expected that cloud mathematical services will lead to a fundamental change in the entire education system and a change in the status of mathematical knowledge in modern society.

Mathematical knowledge can be used effectively, while avoiding complex calculations and many typos and errors.
The MathPartner service appeared in 2011 and is the first cloud mathematics system. Today this service is available at  \cite{MP1}. 

You can use any modern browser to access the server.

The main object on the server is the user's notebook. Texts in a notebook can be loaded and unloaded as regular text files or edited as in a regular editor. 

 This article describes MathPartner and the new capabilities it provides for engineers and scientists, math and science teachers, and their students.

We highlight three generations of symbolic-numeric software and explore how MathPartner differs from other cloud-based tools of the latest generation. We explain in detail its advantages for education at school and at university, those features that make education more intensive and students’ knowledge deeper.
 
In this short article it is impossible to give any complete description of cloud mathematics MathPartner. 
 
Therefore, for a more complete introduction to this service, we can recommend other articles \cite{M17}-\cite{MS21A}, as well as the user manual \cite{UM} and a list of publications  before 2015 \cite{n0} that were used in the project. In addition, the service is equipped with a detailed “Help” section, where the user manual is located in an interactive form, where you can familiarize yourself with the capabilities of the service by running training examples. Alternatively, you can simply copy these examples and transfer them into your workbook. 
 
In the second section, we look at what changes have occurred in the teaching of mathematics and natural sciences in recent years.

The third section is devoted to a description of the Mathpar language.

The fourth section is devoted to a description of the library of mathematical texts in the Mathpar language.

The fifth section describes the technology for preparing and delivering lectures.

 The sixth section describes the process of student mastering new material.

 In section seven, we take a retrospective look at the history of symbolic numerical computing. We distinguish three generations of computer algebra systems.

The eighth section is devoted to the differences between MathPartner cloud mathematics and other third-generation systems.

In the ninth section we consider one problem from a school physics course.
In this example, we show what intensifies the educational process and in-depth study of the subject.

The tenth section summarizes the results.

\section{What's new in natural science educational technologies?}

 In the humanities, information technology has made a genuine revolution. Essentially, humanities lectures will die out. They will be replaced by videos that will be recorded by the best lecturers who will create series of lectures on the humanities. Practical classes can include live dialogues, automated tests, and educational essays that can be checked by artificial intelligence.

 What's new today in the teaching of mathematics and natural sciences, that is, those disciplines in which mathematics acts both as a language and as a tool?

 The surprising fact is that today there is still no generally accepted way of storing mathematical text that can not only be reproduced visually, but could be unambiguously interpreted and “computed.” Textbooks presenting the theory exist separately, just as they did 200 years ago, and the corresponding practice exists separately and uses a notebook and pen.

 Today, almost every developed text editor allows you to insert mathematical expressions into text. There is also a special editor for obtaining high-quality mathematical text - this is the TeX editor, which was created by Donald Knuth almost 40 years ago. Today, editors created on its basis are popular, such as AMSTeX, LaTex and others more or less complete. Their important feature is that mathematical expressions are highlighted with special brackets and, thus, clearly separated from ordinary text. In addition, they offer two playback methods. One method is intended for editing, and the other is intended for traditional reproduction of mathematical text.

 In editing mode, you can change the text, and you can see special characters designed to depict mathematical expressions and you can see special brackets that separate mathematical text.

 In the traditional mode of reproduction of mathematical text, all these auxiliary special characters are hidden, and the mathematical expression is reproduced in a standard form.
 For example, this article was prepared in the LaTeX editor. In order to reproduce integral $\int f(x) d x $ in the text, you need to insert the following fragment: \$$\backslash$int f(x)dx\$.

 A lecturer can create a high-quality lecture by presenting traditional mathematical notation, but that's about it.Such text remains “dead” text. You cannot use previously introduced mathematical objects, you cannot calculate mathematical expressions, or give examples of calculations.

 And of course, today it is easier for a student to solve problems by writing down the solution on paper than to type and save the solution, like an article in a scientific journal.

Thus, we see that the process of mastering and applying natural science disciplines has changed very little.
  
 \section{Language Mathpar }
 
 Mathpar is a subset of the LaTeX language that is used in MathPartner cloud mathematics.

At the same time, its main part is intended for storing mathematical text, which is not only reproduced visually, but is clearly interpreted and calculated.
\subsection{Example}
If the entry looks like this: \\
\centerline{\bf  a = 2; f = a $\backslash$cos(2x); g=$\backslash$int(f) d x;}, \\
 then in addition to the fact that it will be reproduced in this form \\
$$a=2; f=a\cos(2x); g=\int(f) d x, $$
  all expressions will be evaluated and the calculated value of the last expression will appear in the output field: 
  $$\sin(2x).$$
\subsection{Example}
 If you need to find out how several expressions were calculated, then you need to use the $\backslash  print()$ operator. For example, the operator $\backslash print(f,g)$ will cause when you type \\
\centerline{\bf a=2; f=a$\backslash$cos(2x); g= $\backslash$int(f) d; $\backslash$print(f,g);}
the following will be output:
$$f= 2\cos(2x) ; \ g=\sin(2x) .$$
All entered symbols and their meanings will be saved in the user's work page and can be used by him as long as his session remains open.

Unlike LaTeX, Mathpar distinguishes between active and passive parts of the text. LaTeX has only a passive part, 
it consists of text and mathematical expressions, which are separated by special brackets. The dollar symbol is usually 
used as such parentheses.

In Mathpar, the passive part is enclosed in double quotes. 
Anything not enclosed in double 
quotes is active text. These are the expressions that need to 
be calculated, and their values must be saved and used in the
future. They are separated from each other by a semicolon or other passive text in double quotes. 

In the passive part of the text, as in LaTeX, you can find mathematical expressions, but here they should be highlighted with dollar signs.
\subsection{Example with passive and active parts}
"We can first evaluate the indefinite integral,"

$a=2; f=a\backslash cos(2x); g=\backslash int(f) d x;$

"and then check that we haven't made a mistake "

$h=  \backslash D\_\{x\}(g); ; \backslash print(g,h) ;$

OUT: 

$g=\backslash sin(2x) $

 $h= 2\backslash cos(2x) ; $
 
 In image mode it looks like this:

OUT: 

$g=\sin(2x) $

 $h= 2\cos(2x) ; $

 Thus, any reasonable mathematical text that is written in the Mathpar language can not only be visualized, but all its expressions in the active part can be calculated.

This seemingly small addition is, in fact, the biggest breakthrough for all areas of human activity where mathematics is used.

Any fragment of Mathpar text that you need can simply be copied and pasted into the browser window that the user opened on the MathPartner service. And this fragment will be immediately calculated.

  \section{Library of mathematical texts }

The GitHub web service contains a library of texts in the Mathpar language. It is available at http://github.com/mathpar.
You can view text files directly on the site, but you can also download the entire archive.

The service supports receiving and editing code via Git and SVN. You can register with the GitHub web service and then post your own texts or receive texts posted by others. You can also download the archive of the entire library as a zip archive: https://github.com/mathpar/mathpar/archive/master.zip.

Archive root address\\
https://github.com/mathpar/mathpar/tree/master/LIBRARY/mathpar.

For example, if you choose this path EUROPE/United$\_$Kingdom-uk/ \\
SCHOOL/MATH/9-grade/Book-A001/ENSM09A0010517.txt
 then you will receive a section of the school textbook on mathematics
 "Trigonometric functions, numeric circle."
 
 If you choose this path: \\
EUROPE/Ukraine-ua/UNIVERSITY/ MATH/Computer$\_$Algebra/ \\
Malaschonok/Computer$\_$Algebra$\_$Chap$\_$01/101$\_$lect$\_$CA.txt,
 then all paragraphs of the first section of the computer algebra textbook will be available.

 These files can be saved and then loaded into the workbook window on the MathPartner service.

  \section{ Preparation and delivery of lectures}

The teacher can accompany the educational material with examples. He transfers files with lectures to students. Students can copy these examples from lectures. They can use them for practical mastery of new material and for independently solving new problems.

Let's look at how a lecturer can organize lectures. The classroom where the lecture is taking place must have a computer connected to a projector and with Internet access. Using a browser, the lecturer should go to service MathPartner 
and go to the user workbook.
 
In the left sidebar, the lecturer can open the “Files” menu and select the top line “Download text”. The path to the files on his computer will open to him. He must indicate the file with the text of the lecture. This lecture file will be loaded as a sequence of separate windows.

 Above each window on the left there are three control signs:
 
“triangle” – start execution,

 “switch” – switch image mode to text mode,
 
“plus” – add a new empty window below.

Launching for execution occurs for each window separately. There are two possibilities for this. You can either click on the start button on the screen or type the key combination Cntr-Enter.
After launching execution, the result is output, and the text mode changes to image mode. To return to the text editing mode, just use the switch or simply click the mouse in the desired window.

The lecturer can sequentially show fragments of the lecture, divided into separate windows, by sequentially pressing the “execute” buttons and commenting on the image on the screen. He can change ready-made examples to new ones and get different results. Thus, the student learns to solve such problems.

The lecturer has the opportunity to accompany the lecture by constructing drawings, graphs, depicting surfaces and other graphic objects. Such operators are part of the language Mathpar. It is possible to construct a series of graphs and demonstrate the “small movie” that is obtained from these graphs, as from individual frames.

It is convenient to use such lectures for remote learning. To do this, the lecturer must open an image of the lecture on his personal computer and show his screen.

  \section{Mastering new material by the student}

The student can receive all files with lecture texts. After listening to a lecture at the university, he can run all the examples that were demonstrated in the lecture at home. To master new material, he must learn to solve problems. It is impossible for a student to master new material without independently solving problems. Therefore, you need to have a problem book from which you can download the texts of educational problems in the Mathpar language. To solve problems, he can use the MathPartner service, just as the lecturer did.

 He can easily copy the necessary operators from lectures. The student can avoid performing routine arithmetic operations and operations with cumbersome expressions due to the functions built into MathPartner.

Student can save all texts for solving problems and their results as a text file. To do this, in the left sidebar he can open the “Files” menu and select the second line “Save text”. The text is saved in the download folder.

If the student was unable to solve the problem, then he can send his file by mail to the teacher and ask for help and clarification. The teacher can upload a file to the service and detect student errors, using the help of the MathPartner service. In this case, the teacher is freed from routine actions when checking students’ work.

Of course, control work can also be carried out in this mode. Having completed the test, the student sends a file that contains complete solutions to the problems to the teacher.

  \section{ Hindsight to the symbolic-numerical calculations}

Let's look at how Mathpartner came into being, what systems preceded it, and how it fundamentally differs from them.

The field of computer science known as symbolic numerical computing is also called analytical computing or computer algebra.

Three generations can be distinguished in the history of the development of computer algebra systems.

\subsection{The first generation}
The first generation of computer algebra systems appeared in the 1960s.

The very first systems of computer algebra are Schoonschip and Analitik.

 Schoonschip was developed by Martinus J. G. Veltman while working at the Stanford Laboratory. The initial version, dating to December 1963. This system was used by his student Gerardus 't Hooft. He proved the possibility of renormalization of non-Abelian gauge theories. In 1999, M. Veltman and G. ’t Hooft were awarded the Nobel Prize in Physics.\cite{n1}

 The Analitik language for MIR-1 computers was developed in 1964 in Kyiv under the leadership of Viktor Glushkov \cite{n2}.

\subsection{The second generation}
The second generation of computer algebra systems are graphical systems, which appeared in the early 1990s. This generation is associated with the advent of commercial systems Mathematica and Maple, which had a graphical shell. They remain the most common tools for users of computer algebra systems today. A little bit  later, symbolic calculations became available to users of the MATLAB system.

These systems significantly expanded the areas of application of symbolic computing. Many teachers and students of technical and natural sciences, many engineers and scientists began to use them.

 At the same time, the share of users of computer algebra systems in the total number of those who study and apply mathematics remains insignificant. The main reasons are difficulty in mastering and isolation from the school and university educational process.

\subsection{The thired generation}
 The third generation is cloud-based tools for symbolic-numeric computing. They appear in the 2010s.

 These include Wolfram Alpha, MathPartner, Symbolab, and other lesser known ones.

Let's try to figure out how MathPartner differs from other third-generation systems and why we believe that this is a breakthrough technology.
 
 WolframAlpha is an answering system developed by Wolfram Research \cite{n3}.   It is offered as an online service that answers actual queries by calculating answers based on data from external sources.

WolframAlpha was released in 2009 and is based on the  Wolfram Mathematica. WolframAlpha collects data from academic and commercial websites such as the CIA World Factbook, USGS, and others.

 Symbolab is a system that is designed to solve mathematical problems in a number of subjects. Type of site: Answer engine \cite{n4}. The user receives the answer to the problem and can get a step-by-step solution. It was developed by the Israeli company EqsQuest Ltd. in 2011. In 2020, the company was acquired by US edtech website Course Hero \cite{n5}.

Thus, we see that these services are designed to obtain answers to questions that relate to some branches of mathematics.

This is a very convenient tool for a student who does not want to solve his homework himself, but wants to find a ready-made one solution. In some cases, he may even receive a ready-made step-by-step solution.

How does this affect the quality of education? What is new in the educational process?

You can remember that when there were no computers and the Internet, problem books were used, which contained detailed solutions to all problems.
The student could rewrite the solution from such a problem book. Now we essentially have the same problem book, only of a very large size. And it is not in paper, but in electronic form.
One student is trained on a small number of standard problems. The fact that there are many more problems that have ready-made solutions has little effect on the quality of training.

Let's see what's new with MathPartner. Why do we differentiate it from other third generation systems?

\section{What's new with MathPartner?}

\subsection{ This is a new form for storing mathematical knowledge}

 Mathpar is a new form of mathematical writing. Previously, mathematical text could be reproduced on a computer using various editors. The most popular are MS Office Writer and LaTeX. Now all users can write in Mathpar and take advantage of LaTeX. You can reproduce any mathematical texts on it. At the same time, a new quality appeared. Mathematical operators can either be simply displayed or automatically calculated where required.

\subsection{This is an effective way to apply mathematical knowledge }
 
We are accustomed to having an arithmetic calculator on both our phones and computers.

Wolfram Alpha allows you to quickly evaluate a single operator online that Wolfram Mathematica can evaluate.

MathPartner allows you to solve complex multi-step problems online, save and then use intermediate results.

But now the solution to any problem is at the same time a traditional mathematical text that is easy to read, save, forward, and make publicly available.

Anyone can use this solution. He will be able, if necessary, to change this solution and use it to solve another problem.

 In addition, the MathPartner service tools provide many hints and eliminate many routine calculations.

 Any person who has some mathematical knowledge will be able to effectively and quickly apply their knowledge in practice

\subsection{ It is an effective educational technology for learning mathematics and natural science  subjects}

 The process of teaching mathematics, physics and technical disciplines is radically intensified.

Really. If we study a new section of mathematics, then we rely on the knowledge of previous sections, which we already master both theoretically and practically, using a certain set of operators.

When studying a new section, we will freely use this set of operators, which will allow us to quickly solve problems from the new section. We will now be able to solve many more problems, which means we will study the new section more deeply.

Let's look at an example.

\section{ Example of solving a problem in physics}

Let's consider an example of solving physical problems on the topic "Heat Transfer".

\subsection{What's new in the formulation of the educational task}

The student opens the electronic textbook and sees the following text of the problem conditions.\\
\\
"EXERCISE 1 \\
A piece of ice which has mass"\\
M = 10 kg\\
"was put in a vessel. The ice has temperature "\\
T = -10 $ \backslash degreeC$ ;\\
"Find the mass of water in a vessel after transferring the"\\
q = 20000 kJ\\
"amount of heat. Specific heat of water heating is equal "\\
c\_v = 4.2 kJ/(kg $; \backslash degreeC$) ; \\
"Specific heat of ice heating is equal"\\
c\_i = 2.1 kJ/(kg $ \backslash degreeC$ );\\
"The heat of fuson of ice is equal" \\
r = 330 kJ/kg; \\
"Specific heat of vaporization of water is equal" \\
 $\backslash lambda$ = 2300 kJ/kg; \\
\\
The student starts the “execution”, after which he sees such text.\\
\\
EXERCISE 1 \\
A piece of ice which has mass\\
$M = 10 kg$\\
was put in a vessel. The ice has temperature \\
$T = -10 C^{o}$;\\
Find the mass of water in a vessel after transferring the\\
$q = 20000 kJ $\\
amount of heat. Specific heat of water heating is equal\\
$c_v = 4.2 kJ/(kg\ C^{o})$;\\
Specific heat of ice heating is equal\\
$c_i = 2.1 kJ/(kg C^{o}) $;\\
The heat of fusion of ice is equal\\
$r = 330 kJ/kg $; \\
Specific heat of vaporization of water is equal \\
$\lambda = 2300 kJ/kg; $\\
\\
In this case, all notations are automatically entered into the runtime environment.

\subsection{What's new in solving a learning problem}

 The student does not have to enter this data again. He needs to show knowledge of physical laws, write down the necessary equations and solve them. This is how he writes down the solution and explains it in the comments.\\
\\
"SOLUTION OF EX. 1. \\ 
\\
q\_1 = M c\_i (0 - T) " The heat for heating ice"  \\
q\_2 = M r " The heat to melt ice"  \\
q\_3 = M c\_v(100 $\backslash$degreeC) " The  heat for heating water"\\
"Let x be the unknown mass of the remaining water " \\
 q\_4 = (M - x) $\backslash$lambda " The  heat to evaporate water" \\
" We have to solve the equation with  unknown mass x: " \\
mass = $\backslash$solve(q = q\_1 + q\_2 + q\_3 + q\_4); \\
$\backslash$print(mass);\\
\\
The $\bf solve$ operator solves a linear algebraic equation with unknown $x$. We used the fact that by default variable $x$  is considered unknown. It is clear that we can change the settings if desired.

The units of physical quantities act as ordinary literal variables and are automatically reduced where possible.

The student executes his solution, after which the text takes on this form.\\
\\
SOLUTION OF EX. 1. \\
\\
$q_1 = M c_i (0 - T)$  The heat for heating ice \\
$q_2 = M r$  The heat to melt ice\\
$q_3 = M c_v (100 C^{o})$  The  heat for heating water\\
Let x be the unknown mass of the remaining water  \\
$q_4 = (M - x)\lambda$  The  heat to evaporate water \\
 We have to solve the equation with  unknown mass x:   \\
$mass = {\bf solve}(q = q_1 + q_2 + q_3 + q_4); \\
{\bf print}(mass);\\
 out:\\
mass = 10710 \cdot kg/ 2300 
$

If a student wants to simplify this fractional rational expression, then he can enter the $value$ operator in the next window and get an approximate solution: \\
$mass = \backslash value(mass);\\
 out: \\
 4.66 \cdot kg \\
$  
\section{Conclusion }
 
 The era of “manual labor” in mathematics ended with the advent of the Mathpar language, which provides automatic calculation of mathematical expressions in the MathPartner service.

Such automation of mathematics will give a powerful impetus to both mathematics and natural science education, and to those areas of science and technology where mathematics is actively used.

Cloud mathematics MathPartner - new
 breakthrough technology. This is undoubtedly a very important tool that paves the way from the modern information society to the knowledge society.   
\begin{credits}
\subsubsection{\ackname}

 The first author was supported by a grant from the Simons Foundation, ID: 1290592 (2024).

Cloud mathematics MathPartner is hosted on the server of the National University of Kiev-Mohyla Academy.

 The first author expresses gratitude to the Mathematica Visiting Scholar Grant Program (2000-2001) and  USA Civilian Research and Development Foundation, TGP-352 (2000). They served as the impetus for the start of the MathPartner project.

 \subsubsection{\discintname}
 The authors have no competing interests to declare that are
relevant to the content of this article. 
 \end{credits}
%
%
%
%

\end{document}